\documentclass[preprint,12pt]{elsarticle}

\usepackage{amssymb}
\usepackage{amsthm}
\usepackage{graphics}
\usepackage{epsfig}
\usepackage{graphicx}
\usepackage{color}
\usepackage{hyperref}
\usepackage{a4wide}
\usepackage{amsmath}
\usepackage{amsfonts}
\usepackage{enumerate}
\newcommand{\pf}{\noindent {\bf Proof: }}

\newtheorem{lemma}{\bf Lemma}[section]
 
\newtheorem{theorem}{\bf Theorem}[section]

\newtheorem{corollary}[lemma]{\bf Corollary}

\newtheorem{remark}{\bf Remark}[section]

\journal{Journal of  Algebra and its Applications}

\begin{document}

\begin{frontmatter}



\title{On Non-Zero Component Graph of Vector Spaces over Finite Fields}



\author{Angsuman Das\corref{cor1}}
\ead{angsumandas@sxccal.edu}

\address{Department of Mathematics,\\ St. Xavier's College, Kolkata, India.\\angsumandas@sxccal.edu}
\cortext[cor1]{Corresponding author}


\begin{abstract}
In this paper, we study non-zero component graph $\Gamma(\mathbb{V})$ on a finite dimensional vector space $\mathbb{V}$ over a finite field $\mathbb{F}$. We show that the graph is Hamiltonian and not Eulerian. We also characterize the maximal cliques in $\Gamma(\mathbb{V})$ and show that there exists two classes of maximal cliques in $\Gamma(\mathbb{V})$. We also find the exact clique number of $\Gamma(\mathbb{V})$ for some particular cases. Moreover, we provide some results on size, edge-connectivity and chromatic number of $\Gamma(\mathbb{V})$.
\end{abstract}


\end{frontmatter}


\section{Introduction}
The study of graphs associated with of various algebraic structures was initiated by Beck \cite{beck} who introduced the idea of zero divisor graph of a commutative ring with unity.  Till then, a lot of research, e.g., \cite{survey2,zero-divisor-survey,anderson-livingston,graph-ideal,power1,power2,mks-ideal,int-vecsp-2,int-vecsp-1} has been done in connecting graph structures to various algebraic objects like semigroups, groups, rings, vector spaces etc. In this paper, we continue the study of one such graph called Non-Zero Component Graph \cite{angsu-comm-alg} of a finite dimensional vector space $\mathbb{V}$ over a field $\mathbb{F}$ with respect to a basis $\{\alpha_1,\alpha_2,\ldots,\alpha_n\}$ of $\mathbb{V}$.

Throughout this paper, $\mathbb{F}$ is a finite field with $q=p^t$ elements, $p$ being a prime and $n=dim_{\mathbb{F}}(\mathbb{V})$. We study the structures of size, connectivity, maximal cliques, clique number, chromatic number and Hamiltonicity of these graphs and other related concepts.

\section{Some Preliminaries}
In this section, for convenience of the reader and also for later use, we recall some definitions, notations and results concerning elementary graph theory. For undefined terms and concepts the reader is referred to \cite{west-graph-book}.

By a graph $G=(V,E)$, we mean a non-empty set $V$ and a symmetric binary relation (possibly empty) $E$ on $V$. The set $V$ is called the set of vertices and $E$ is called the set of edges of $G$. Two element $u$ and $v$ in $V$ are said to be adjacent if $(u,v) \in E$. $H=(W,F)$ is called a {\it subgraph} of $G$ if $H$ itself is a graph and $\phi \neq W \subseteq V$ and $F \subseteq E$. If $V$ is finite, the graph $G$ is said to be finite, otherwise it is infinite. If all the vertices of $G$ are pairwise adjacent, then $G$ is said to be {\it complete}. A complete subgraph of a graph $G$ is called a {\it clique}. A {\it maximal clique} is a clique which is maximal with respect to inclusion. The {\it clique number} of $G$, written as $\omega(G)$, is the maximum size of a clique in $G$. The {\it chromatic number} of $G$, denoted as $\chi(G)$, is the minimum number of colours needed to label the vertices so that the adjacent vertices receive different colours. A subset $I$ of $V$ is said to be {\it independent} if any two vertices in that subset are pairwise non-adjacent. The {\it independence number} of a graph, $\alpha(G)$, is the maximum size of an independent set of vertices in $G$. A {\it path} of length $k$ in a graph is an alternating sequence of vertices and edges, $v_0,e_0,v_1,e_1,v_2,\ldots, v_{k-1},e_{k-1},v_k$, where $v_i$'s are distinct (except possibly the first and last vertices) and $e_i$ is the edge joining $v_i$ and $v_{i+1}$. We call this a path joining $v_0$ and $v_{k}$. A graph is {\it connected} if for any pair of vertices $u,v \in V,$ there exists a path joining $u$ and $v$. A path with first and last vertices same is called a {\it cycle}. A graph is said to be {\it Hamiltonian} if it contains a cycle consists of all the vertices in $G$. A graph is said to be {\it Eulerian} if it contains a cycle consists of all the edges in $G$ exactly once.
\section{Definitions and Some Basic Results}
Firstly, we recall the definition of Non-zero Component graph of a finite dimensional vector space and some preliminary results from \cite{angsu-comm-alg}.

\begin{figure}[ht]
\centering
\begin{center}

$\begin{array}{lr}
\includegraphics[scale=.4]{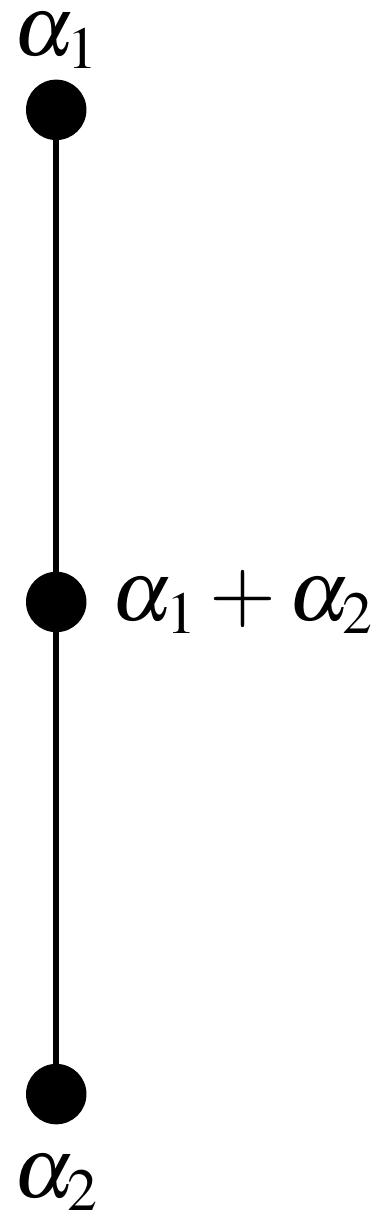}~~~~~~~~~~~ & ~~~~~~~~~~ \includegraphics[scale=.4]{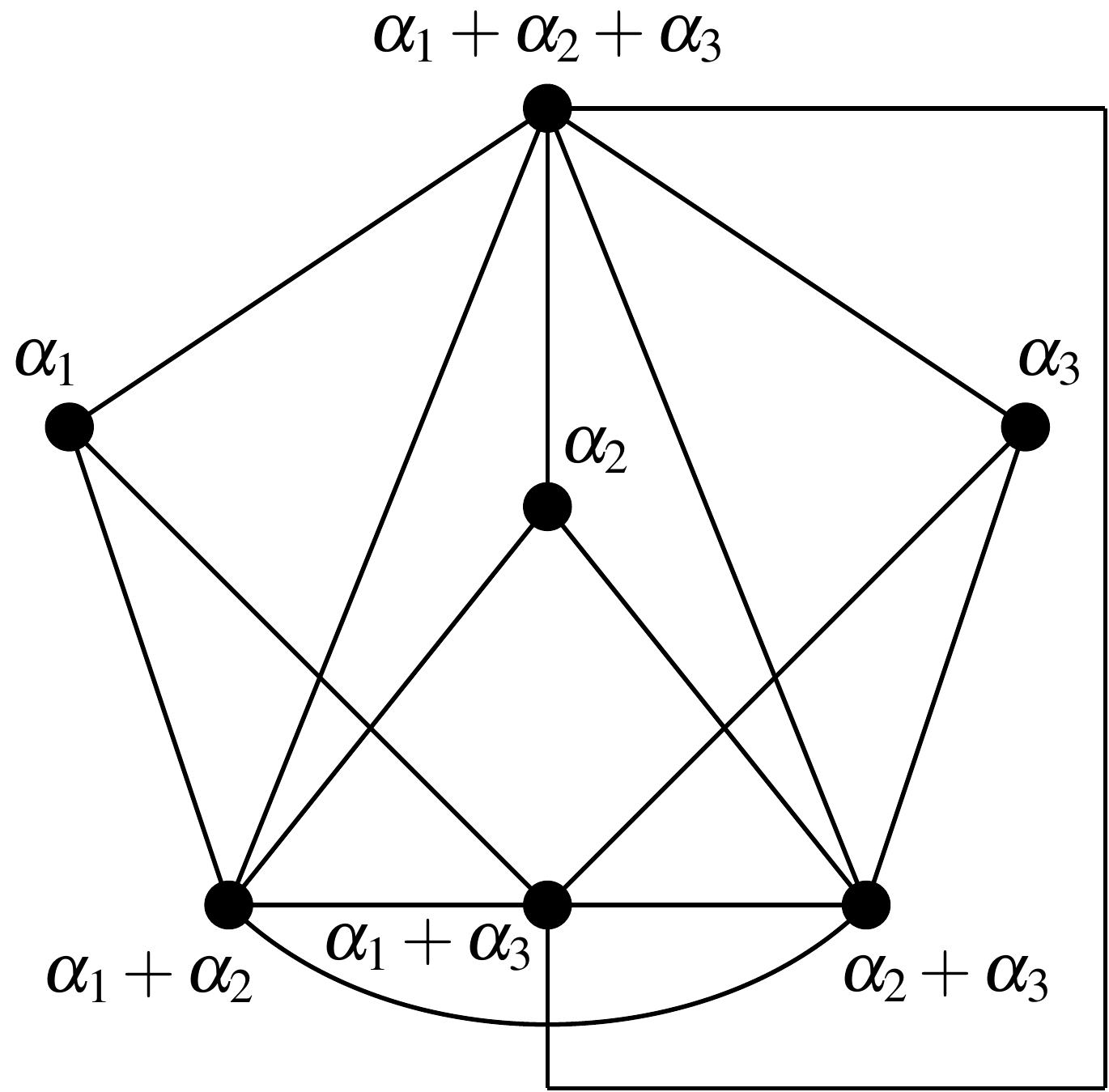}\\

\mathbf{dim(\mathbb{V})=2; \mathbb{F}=\mathbb{Z}_2} & \mathbf{dim(\mathbb{V})=3; \mathbb{F}=\mathbb{Z}_2}
\end{array}
$

\vspace{.3in}
$
\begin{array}{c}
\includegraphics[scale=.4]{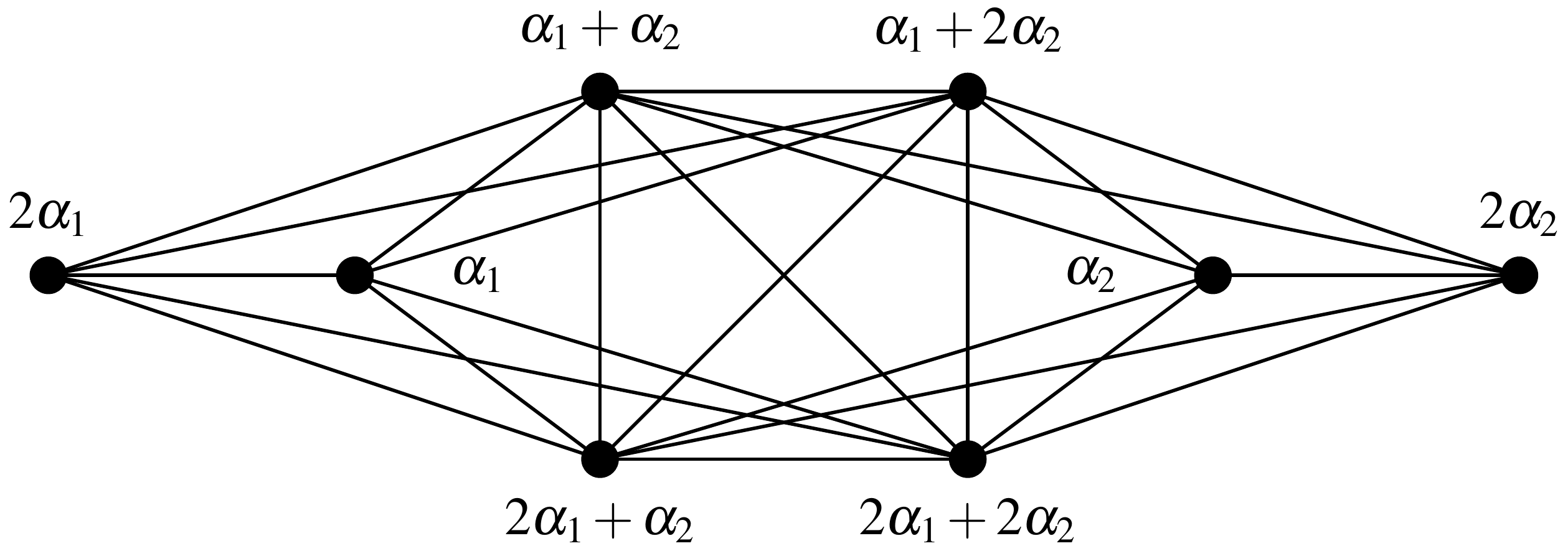}\\

\mathbf{dim(\mathbb{V})=2; \mathbb{F}=\mathbb{Z}_3}
\end{array}
$
\caption{Examples of $\Gamma(\mathbb{V})$}
\label{example-figure}
\end{center}
\end{figure}

Let $\mathbb{V}$ be a vector space over a field $\mathbb{F}$ with $\{\alpha_1,\alpha_2,\ldots,\alpha_n\}$ as a basis and $\theta$ as the null vector. Then any vector $\mathbf{a} \in \mathbb{V}$ can be expressed uniquely as a linear combination of the form $\mathbf{a}=a_1\alpha_1+a_2\alpha_2+\cdots+a_n\alpha_n$. We denote this representation of $\mathbf{a}$ as its basic representation w.r.t. $\{\alpha_1,\alpha_2,\ldots,\alpha_n\}$. We define {\it Non-Zero Component graph} of a finite dimensional vector space $\Gamma(\mathbb{V}_\alpha)=(V,E)$ (or simply $\Gamma(\mathbb{V})$) with respect to $\{\alpha_1,\alpha_2,\ldots,\alpha_n\}$ as follows: $V=\mathbb{V}\setminus \{\theta\}$ and for $\mathbf{a},\mathbf{b} \in V$, $\mathbf{a} \sim \mathbf{b}$ or $(\mathbf{a},\mathbf{b}) \in E$ if $\mathbf{a}$ and $\mathbf{b}$ shares atleast one $\alpha_i$ with non-zero coefficient in their basic representation, i.e., there exists atleast one $\alpha_i$ along which both $\mathbf{a}$ and $\mathbf{b}$ have non-zero components. Unless otherwise mentioned, we take the basis on which the graph is constructed as $\{\alpha_1,\alpha_2,\ldots,\alpha_n\}$. For some examples of $\Gamma(\mathbb{V})$, see Figure \ref{example-figure}.

{\theorem \cite{angsu-comm-alg} \label{diameter-theorem} $\Gamma(\mathbb{V})$ is connected and $diam(\Gamma)=2$.}
{\theorem \cite{angsu-comm-alg} $\Gamma(\mathbb{V})$ is complete if and only if $\mathbb{V}$ is one-dimensional.}
{\theorem \cite{angsu-comm-alg} \label{independence-number-theorem} The independence number of $\Gamma(\mathbb{V})$, $\alpha(\Gamma(\mathbb{V}))=dim(\mathbb{V})$.}

{\theorem \label{degree-theorem} \cite{angsu-comm-alg} Let $\mathbb{V}$ be a vector space over a finite field $\mathbb{F}$ with $q$ elements and $\Gamma$ be its associated graph with respect to a basis $\{\alpha_1,\alpha_2,\ldots,\alpha_n\}$. Then, the degree of the vertex $c_1\alpha_{i_1}+c_2 \alpha_{i_2}+\cdots+c_k \alpha_{i_k}$, where $c_1 c_2 \cdots c_k \neq 0$, is $(q^k -1)q^{n-k}-1$.}

 Now we prove some basic results when the underlying field $\mathbb{F}$ is finite.

{\theorem $\Gamma$ is not Eulerian.}\\
\pf  If $q$ is an odd prime then $\Gamma$ is not Eulerian as from Theorem \ref{degree-theorem} every vertex is of odd degree. If $q=2$, then from Theorem \ref{degree-theorem}, all the vertices with $1<k<n$ are of odd degree. Thus the graph is not Eulerian in any case.\qed

{\lemma \label{min-degree-lemma} If $\mathbb{V}$ be a $n$-dimensional vector space over a finite field $\mathbb{F}$ with $q$ elements, then the minimum degree $\delta$ of $\Gamma(\mathbb{V})$ is $q^{n-1}(q-1)-1$.}\\
\pf From Theorem \ref{degree-theorem}, the degree of the vertex $c_1\alpha_{i_1}+c_2 \alpha_{i_2}+\cdots+c_k \alpha_{i_k}$, where $c_1 c_2 \cdots c_k \neq 0$, is $(q^k -1)q^{n-k}-1$ i.e., $q^n-q^{n-k}-1$. Thus the degree will be minimized if $k=1$ and hence $\delta=q^n-q^{n-1}-1$.\qed

{\corollary \label{edge-connectivity} Edge connectivity of $\Gamma(\mathbb{V})$ is $q^{n-1}(q-1)-1$.}\\
\pf From \cite{plesnik}, as $\Gamma(\mathbb{V})$ is of diameter $2$ (by Theorem \ref{diameter-theorem}), its edge connectivity is equal to its minimum degree, i.e., $q^{n-1}(q-1)-1$. \qed

{\theorem \label{size-theorem} If $\mathbb{V}$ be a $n$-dimensional vector space over a finite field $\mathbb{F}$ with $q$ elements, then the order of $\Gamma(\mathbb{V})$ is $q^n-1$ and the size $m$ of $\Gamma(\mathbb{V})$ is $$\dfrac{q^{2n}-q^n+1-(2q-1)^n}{2}.$$}\\
\pf It is trivial to observe that the order of $\Gamma(\mathbb{V})$ is $q^n-1$. By Theorem \ref{degree-theorem}, the degree of the vertex $c_1\alpha_{i_1}+c_2 \alpha_{i_2}+\cdots+c_k \alpha_{i_k}$, where $c_1 c_2 \cdots c_k \neq 0$, is $(q^k -1)q^{n-k}-1$. Now, there are $\binom{n}{k}(q-1)^k$ vectors with exactly $k$ many $\alpha_i$'s in its basic representation. Since, $2m=$ sum of the degree of all the vertices in $\Gamma(\mathbb{V})$, we have
$$2m=\sum_{k=1}^n \binom{n}{k}(q-1)^k\left[(q^k -1)q^{n-k}-1\right]=\sum_{k=1}^n \binom{n}{k}(q-1)^k\left[(q^n -1)-q^{n-k}\right]$$
$$=(q^n-1)\sum_{k=1}^n \binom{n}{k}(q-1)^k - \sum_{k=1}^n \binom{n}{k}(q-1)^k q^{n-k}$$
$$=(q^n-1)\left[(q-1+1)^n-1 \right]-\left[(q+q-1)^n - q^n \right]=(q^n-1)^2 + q^n - (2q-1)^n$$
$$=q^{2n}-q^n+1-(2q-1)^n$$ and hence the result. \qed

\section{Maximal Cliques in $\Gamma(\mathbb{V})$}
In this section, we study the structure of maximal cliques in $\Gamma(\mathbb{V})$ and find its clique number. In the process, we show that $\Gamma(\mathbb{V})$ possess two different classes of maximal cliques.

Let $S_\beta$ (skeleton of $\beta$) be the set of $\alpha_i$'s with non-zero coefficients in the basic representation of $\beta$ with respect to $\{\alpha_1,\alpha_2,\ldots, \alpha_n\}$. It is to be noted that two distinct $\beta$ may have same $S_\beta$. Also $1 \leq |S_\beta| \leq n, \forall \beta \in \Gamma(\mathbb{V}_\alpha)$. Let $M$ be a maximal clique in $\Gamma(\mathbb{V}_\alpha)$ and define $S(M)=\{S_\beta: \beta \in M\}$ and $S[M]=\{|S_\beta|: S_\beta \in S(M)\}$. Since $M$ is a clique, $S_\alpha \cap S_\beta \neq \emptyset, \forall \alpha, \beta \in M$. By maximality of $M$, if $\alpha \in M$ and $S_\alpha \subset S_\beta$ for some $\beta \in \Gamma(\mathbb{V}_\alpha)$, then $\beta \in M$. 

As $\emptyset \neq S[M] \subset \mathbb{N}$, by well-ordering principle, it has a least element, say $k$. Then there exist some $\beta^* \in M$ with $|S_{\beta^*}|=k$, where $\beta^*=c_1\alpha_{i_1}+c_2\alpha_{i_2}+\cdots+c_k\alpha_{i_k}$. Now accordingly as $k \leq n/2$ or $k>n/2$, we show that there exists two types of maximal cliques in $\Gamma(\mathbb{V}_\alpha)$.

{\theorem \label{k-leq-n/2-clique} Let $M$ be a maximal clique in $\Gamma(\mathbb{V}_\alpha)$. If $k$ is the least element of $S[M]$ and $k \leq n/2$, then $M$ belongs to a family of maximal cliques $\{M_{k,i}: 1 \leq k \leq n/2 ;i \in \{1,2,\ldots , n\}\}$ of $\Gamma(\mathbb{V}_\alpha)$ where $M_{k,i}=\{\beta \in \Gamma(\mathbb{V}_\alpha): \alpha_i \in S_\beta \mbox{ and }|S_\beta|\geq k \}$ and $$|M|=(q-1)\sum_{r=k-1}^{n-1} \binom{n-1}{r}(q-1)^r.$$}\\
\pf Since minimum of $|S_\beta|$ for $\beta \in M$ is $k ~(\leq n/2)$ and $S_\alpha \cap S_\beta \neq \emptyset, \forall \alpha, \beta \in M$, by Erdos-Ko-Rado theorem \cite{erdos-ko-rado}, the maximum number of pairwise-intersecting $k$-subsets in $S(M)$ is $\binom{n-1}{k-1}$ and the maximum is achieved only if each $k$-subsets contain a fixed element, say $\alpha_i$. As $M$ is a maximal clique, $M=\{\beta \in \Gamma(\mathbb{V}_\alpha): \alpha_i \in S_\beta \mbox{ and } |S_\beta|\geq k\}$.

Now, the number of $\beta$'s in $M$ with $|S_\beta|=k$ and $\alpha_i \in S_\beta$ is $\binom{n-1}{k-1} (q-1)^k$. Similarly, the numbers of $\beta$'s in $M$ with $|S_\beta|=k+1,k+2,\ldots,n$ and $\alpha_i \in S_\beta$ are $$\binom{n-1}{k} (q-1)^{k+1},\binom{n-1}{k+1} (q-1)^{k+2},\ldots,\binom{n-1}{n-1} (q-1)^n$$ respectively. Thus, we have
$$|M|=\binom{n-1}{k-1} (q-1)^k+ \binom{n-1}{k} (q-1)^{k+1}+ \binom{n-1}{k+1} (q-1)^{k+2}+ \cdots + \binom{n-1}{n-1} (q-1)^n$$
$$=(q-1)\left[ \binom{n-1}{k-1} (q-1)^{k-1}+ \binom{n-1}{k} (q-1)^{k}+ \cdots + \binom{n-1}{n-1} (q-1)^{n-1} \right]$$
$$=(q-1)\sum_{r=k-1}^{n-1} \binom{n-1}{r}(q-1)^r. $$

It is to be noted that for same value of $k$ and by fixing different $\alpha_i$'s, we get different maximal cliques. Since these maximal cliques depends both on $k$ and $\alpha_i$, we get a family of maximal cliques $M_{k,i}$ where $1 \leq k \leq n/2$ and $i \in \{1,2,\ldots , n\}$ and $M \in M_{k,i}$.\qed

{\theorem \label{k>n/2-clique} Let $M$ be a maximal clique in $\Gamma(\mathbb{V}_\alpha)$. If $k$ is the least element of $S[M]$ and $k > n/2$, then $k=\lfloor n/2 \rfloor + 1$ and $M=\{\beta \in \Gamma(\mathbb{V}_\alpha):|S_\beta| \geq \lfloor n/2 \rfloor + 1\}$ and $$|M|=\sum_{r=k}^{n} \binom{n}{r}(q-1)^r.$$}\\
\pf Since $k > n/2$, for all $\alpha, \beta \in M$, $S_\alpha \cap S_\beta \neq \emptyset$. Thus the maximum number of pairwise-intersecting $k$-subsets in $S(M)$ is $\binom{n}{k}$. Thus by similar arguments as in the proof of Theorem \ref{k-leq-n/2-clique},

$$|M|=\binom{n}{k} (q-1)^k+ \binom{n}{k+1} (q-1)^{k+1}+ \cdots + \binom{n}{n} (q-1)^n$$
$$=\sum_{r=k}^{n} \binom{n}{r}(q-1)^r.$$
Now, as $\{\beta \in \Gamma(\mathbb{V}_\alpha):|S_\beta| \geq k + 1\} \subset \{\beta \in \Gamma(\mathbb{V}_\alpha):|S_\beta| \geq k\}$, by maximality of $M$, $M=\{\beta \in \Gamma(\mathbb{V}_\alpha):|S_\beta| \geq k\}$ when $k$ is minimized provided $k>n/2$. Thus $k=\lfloor n/2 \rfloor + 1$ and hence the theorem.\qed

{\remark \label{clique-number-remark} It is obvious from Theorem \ref{k-leq-n/2-clique}, that $|M_{k,i}|$ is maximum when $k=1$, i.e., $M_{1,i}=\{c_1\alpha_1+c_2\alpha_2+\cdots+c_n\alpha_n: c_i \neq 0\}$ and $|M_{1,i}|=(q-1)q^{n-1}$. Thus, the clique number of $\Gamma(\mathbb{V}_\alpha)$ is $$\omega(\Gamma(\mathbb{V}_\alpha))=\max \lbrace (q-1)q^{n-1}, \sum_{r=\lfloor n/2 \rfloor + 1}^{n} \binom{n}{r}(q-1)^r  \rbrace $$ and it depends on the value of $q$ and $n$.
}

{\remark A maximal clique $M$ contains at most one $\alpha_i$, because if $M$ contain $\alpha_i$ and $\alpha_j$, then $\alpha_i \not\sim \alpha_j$ which contradicts that $M$ is a clique. Moreover if $M$ is a maximal clique containing $\alpha_i$, then $M=M_{1,i}$. It follows since every $\beta \in M$ is adjacent to $\alpha_i$, i.e., every $\beta$ has a non-zero component along $\alpha_i$ and hence $M_{1,i} \subset M$. Now, by maximality of $M_{1,i}$, it follows that $M=M_{1,i}$. }

{\corollary \label{clique-number-for-q=2} If $q=2$, the clique number $\omega(\Gamma(\mathbb{V}))=2^{n-1}$.}\\
\pf If $q=2$, $(q-1)q^{n-1}=2^{n-1}$. Now, for $q=2$ and $n$ even (say $2m$), $$\sum_{r=\lfloor n/2 \rfloor + 1}^{n} \binom{n}{r}(q-1)^r=\sum_{r=m + 1}^{2m} \binom{2m}{r}=\dfrac{1}{2}\sum_{r=0}^{2m} \binom{2m}{r} -\binom{2m}{m}<2^{2m-1}=2^{n-1}.$$
If $q=2$ and $n$ is odd (say $2m+1$), $$\sum_{r=\lfloor n/2 \rfloor + 1}^{n} \binom{n}{r}(q-1)^r=\sum_{r=m + 1}^{2m+1} \binom{2m+1}{r}=\dfrac{1}{2}\sum_{r=0}^{2m+1} \binom{2m+1}{r} =2^{2m+1-1}=2^{n-1}.$$
Combining both the cases and using Remark \ref{clique-number-remark}, we get $\omega(\Gamma(\mathbb{V})=2^{n-1}$.\qed

{\corollary \label{clique-number-for-q>2-and-n-odd} If $q>2$ and $n$ is odd, the clique number $$\omega(\Gamma(\mathbb{V}))=\sum_{r=\lfloor n/2 \rfloor + 1}^{n} \binom{n}{r}(q-1)^r.$$}\\
\pf It follows from the inequality proved in Appendix and Remark \ref{clique-number-remark}. \qed

{\corollary \label{chromatic-number-for-q=2} If $q=2$ and $\chi(\Gamma(\mathbb{V})$ be the chromatic number $\Gamma(\mathbb{V})$, then $$2^{n-1} \leq \chi(\Gamma(\mathbb{V}) \leq 2^{n-1}+2^{n-2}-n/2.$$}\\
\pf First part of the inequality follows from Corollary \ref{clique-number-for-q=2} and the fact that for any graph $G$, $\omega(G)\leq \chi(G)$. For the other part, we use the following result from \cite{brigham-dutton}: $$\chi(G)\leq \dfrac{\omega(G)+|G|+1-\alpha(G)}{2}$$ 
where $\alpha(G)$ is the independence number of $G$. Thus, by using Corollary \ref{clique-number-for-q=2} and Theorem \ref{independence-number-theorem}, we get for $q=2$, $$\chi(\Gamma(\mathbb{V})\leq \dfrac{2^{n-1}+2^n-n}{2}=2^{n-1}+2^{n-2}-n/2.$$ \qed

\section{$\Gamma(\mathbb{V})$  is Hamiltonian}
In this section, we prove that $\Gamma(\mathbb{V})$ is Hamiltonian except the case when $q=2$ and $n=2$. Firstly, we recall two classical theorems on sufficient conditions of Hamiltonicity of a graph which are crucial later in our proofs.
{\theorem \label{dirac-theorem} \cite{dirac} {\bf[Dirac]} If $G$ is a connected graph with minimum degree $\delta$ such that $\delta\geq |G|/2$, then $G$ is Hamiltonian.}

{\theorem \label{nash-williams-theorem} \cite{nash-williams} {\bf[Nash-Williams]} If $G$ is a $2$-connected graph with minimum degree $\delta$ and independence number $I$ such that $\delta\geq \max \{(|G|+2)/3,I\}$, then $G$ is Hamiltonian.}

{\theorem If $q>2$, then $\Gamma(\mathbb{V})$ is Hamiltonian.}\\
\pf If $n=1$, $\Gamma(\mathbb{V})$ is complete and hence it is Hamiltonian. So, we assume that $n>1$. Since $q>2$, by principle of mathematical induction it can be shown that for $n\geq 2, q^n> 2q^{n-1}+1$, which implies $2q^n-2q^{n-1}-2>q^n-1$, i.e., $q^n-q^{n-1}-1>(q^n-1)/2$. Now by Lemma \ref{min-degree-lemma} $\delta > |G|/2$ and hence by Theorem \ref{dirac-theorem}, $\Gamma(\mathbb{V})$ is Hamiltonian. \qed

{\lemma \label{2-connected-lemma} If $q=2$ and $n\geq 3$, then $\Gamma(\mathbb{V})$ is $2$-connected.}\\
\pf We prove that removal of any one vertex does not make $\Gamma(\mathbb{V})$ disconnected, i.e., $\Gamma(\mathbb{V}) - \{\alpha\}$ is connected for any non-null vector $\alpha \in \mathbb{V}$. Let $\{\alpha_1,\alpha_2,\ldots,\alpha_n\}$ be a basis of $\mathbb{V}$. Let $\mathbf{a}$ and $\mathbf{b}$ be two arbitrary non-null vectors other than $\alpha$ in $\mathbb{V}$. If they are adjacent in $\Gamma$, we are done. If not, since $\mathbf{a},\mathbf{b} \neq \theta$, $\exists \alpha_i, \alpha_j$ which have non-zero coefficient in the basic representation of $\mathbf{a}$ and $\mathbf{b}$ respectively. Moreover, as $\mathbf{a}$ and $\mathbf{b}$ are not adjacent, $\alpha_i \neq \alpha_j$. Consider $\mathbf{c}=\alpha_i + \alpha_j$. Then, $\mathbf{a}\sim \mathbf{c}$ and $\mathbf{b} \sim \mathbf{c}$ in $\Gamma(\mathbb{V})$. This provides a path between $\mathbf{a}$ and $\mathbf{b}$. However, if $\mathbf{c}$ is the removed vertex $\alpha$, then this path between $\mathbf{a}$ and $\mathbf{b}$ does not exist in $\Gamma(\mathbb{V}) - \{\alpha\}$. But as $n \geq 3$, there exists $\alpha_k$ in the mentioned basis other than $\alpha_i, \alpha_j$. Then consider the vertex $\mathbf{d}=\alpha_i + \alpha_j + \alpha_k$ and observe that $\mathbf{a}\sim \mathbf{d}$ and $\mathbf{b} \sim \mathbf{d}$ in $\Gamma(\mathbb{V}) - \{\alpha\}$. Hence $\Gamma(\mathbb{V}) - \{\alpha\}$ is connected and thereby $\Gamma(\mathbb{V})$ is $2$-connected. \qed

{\lemma \label{delta-inequality-lemma} If $q=2$ and $n\geq 3$, then $\delta\geq \max \{(|\Gamma(\mathbb{V})|+2)/3,\alpha(\Gamma(\mathbb{V}))\}$.}\\
\pf Since $n \geq 3$, $2^{n-1}-4\geq 0 \Rightarrow (2^n -2)+(2^{n-1}-1)\geq 2^n+1 \Rightarrow 3(2^{n-1}-1)\geq 2^n+1$, i.e., $2^{n-1}-1\geq (2^n+1)/3$. Now for $q=2$, $\delta=2^{n-1}-1$ and $|\Gamma(\mathbb{V})|=2^n-1$. Thus, the inequality gives us $$\delta \geq (|\Gamma(\mathbb{V})|+2)/3.$$ Also from Theorem \ref{independence-number-theorem}, $\alpha(\Gamma(\mathbb{V}))=n$ and for $n \geq 3$, $$\delta=2^{n-1}-1 \geq n.$$ Combining the above two inequalities, we get the lemma. \qed

{\theorem If $q=2$ and $n\geq 3$, then $\Gamma(\mathbb{V})$ is Hamiltonian.}\\
\pf The theorem follows from Lemma \ref{2-connected-lemma}, Lemma \ref{delta-inequality-lemma} and Theorem \ref{nash-williams-theorem}. \qed

{\remark For $q=2$ and $n=1$, $\Gamma(\mathbb{V})$ is a single vertex graph and for $q=2,n=2$, $\Gamma(\mathbb{V})$ is isomorphic to $P_3$, a $3$-vertex path (See Figure \ref{example-figure}) and hence not Hamiltonian.}

\section*{Acknowledgement}
The author is thankful to Professor Mridul Kanti Sen for some fruitful suggestions on the paper. The research is partially funded by NBHM Research Project Grant, (Sanction No. 2/48(10)/2013/ NBHM(R.P.)/R\&D II/695), Govt. of India. 


\section*{Appendix}
{\bf A Binomial Inequality:} \label{odd-inequality} If $q>2$ and $n$ is odd, then $$(q-1)q^{n-1}< \sum_{r=\lfloor n/2 \rfloor + 1}^{n} \binom{n}{r}(q-1)^r.$$\\
\pf Let $n=2m+1$. Then $$\sum_{r=\lfloor n/2 \rfloor + 1}^{n} \binom{n}{r}(q-1)^r=\sum_{r=m + 1}^{2m+1} \binom{2m+1}{r}(q-1)^r=(q-1)\sum_{r=m + 1}^{2m+1} \binom{2m+1}{r}(q-1)^{r-1}$$
$$=(q-1)\left[ \binom{2m+1}{m+1}(q-1)^m + \binom{2m+1}{m+2}(q-1)^{m+1} +  \cdots + \binom{2m+1}{2m+1}(q-1)^{2m}  \right]$$
$$>(q-1)\left[ \binom{2m+1}{m}(q-1)^m + \binom{2m+1}{m-1}(q-1)^{m-1} +  \cdots + \binom{2m+1}{0}  \right]$$
~~~~~~~~~~~~~~~~~~~~~~~~~~~~~~~~~~~~~~~~~~~~~~[since $\binom{n}{r}=\binom{n}{n-r}$ and $i < j \Rightarrow (q-1)^i < (q-1)^j$].
$$\mbox{Therefore, }\sum_{r=m + 1}^{2m+1} \binom{2m+1}{r}(q-1)^r > (q-1)\sum_{r=0}^{m} \binom{2m+1}{r}(q-1)^r.~~~~~~~~~~~~~~~~~~~~~~~~~~~~~~~$$
$$\mbox{Now, }q^{2m+1}=(q-1+1)^{2m+1}=\sum_{r=0}^{m} \binom{2m+1}{r}(q-1)^r + \sum_{r=m + 1}^{2m+1} \binom{2m+1}{r}(q-1)^r$$
$$>(q-1+1)\sum_{r=0}^{m} \binom{2m+1}{r}(q-1)^r,$$ 
 $$\mbox{ i.e., }\sum_{r=0}^{m} \binom{2m+1}{r}(q-1)^r < q^{2m}.~~~~~~~~~~~~~~~~~~~~~~~~~~~~~~~~~~~~~~~~~~~~~~~~~~~~~~~~~~~~~~~~~~~~~~~~~~~$$
Thus,
$$\sum_{r=\lfloor n/2 \rfloor + 1}^{n} \binom{n}{r}(q-1)^r=\sum_{r=m + 1}^{2m+1} \binom{2m+1}{r}(q-1)^r=q^{2m+1}-\sum_{r=0}^{m} \binom{2m+1}{r}(q-1)^r$$
$$>q^{2m+1}-q^{2m}=q^{2m}(q-1)=q^{n-1}(q-1).$$\qed

\end{document}